\documentclass[a4paper]{amsart}
\usepackage{amsmath,amssymb,amsthm,amscd} 
\usepackage{graphicx} 

 \newtheorem{thm}{Theorem}[section] 
 
 \newtheorem{prop}[thm]{Proposition}
 \newtheorem{lem}[thm]{Lemma}
\theoremstyle{definition} 
 \newtheorem{defn}[thm]{Definition} 
 \newtheorem*{acknowledgments}{Acknowledgments}
\theoremstyle{remark}  
 \newtheorem{rem}[thm]{Remark}

\newcommand{\cond}[1]{$\mathbf{(#1)}$} 
\newcommand{\A}{{\mathcal A}} 
\newcommand{\B}{{\mathcal B}} 
\newcommand{\C}{{\mathcal C}}

\newcommand{\I}{{\mathcal I}}

\DeclareMathOperator{\Ker}{Ker} 
 
\def\we{\wedge} 
\def\del{\partial} 
 
\newcommand{\card}[1]{ \vert #1 \vert }  
\newcommand{\ideal}[1]{ \langle #1 \rangle }

\begin{document} 

\title[The non-vanishing cohomology]{
The non-vanishing cohomology \\ of Orlik-Solomon algebras} 
\author{Yukihito Kawahara} 
\thanks{Partially supported by Research Fellowships of the Japan Society 
for the Promotion of Science for Young Scientists. }    
\address{Department of Mathematics, Tokyo Metropolitan University \\ 
Minami-Ohsawa 1-1, Hachioji-shi, Tokyo  192-0397, Japan } 
\email{ykawa@comp.metro-u.ac.jp} 
\maketitle 
\begin{abstract} 
The cohomology on the complement of hyperplanes  
with the coefficients in the rank one local system 
associated to a generic weight  
vanishes except in the highest dimension. 
In this paper, 
we construct matroids or arrangements and its weights with  
non-vanishing cohomology of Orlik-Solomon algebras,  
using decomposable relations arising from Latin hypercubes. 
\end{abstract}   

\section{ Introduction } 

Let $R$ be a commutative ring with $1$. 
Write $[n] := \{ 1,2, \ldots, n\}$.  
Let $E = E_R$ denote the graded exterior algebra over $R$ generated 
by $1$ and degree-one elements $e_i$ for $i \in [n]$. 
Define a $R$-linear map $\del : E^p \to E^{p-1}$ by 
$\del{1} =0$, $\del{e_i}=1$ for $i \in [n]$, 
and 
$$ \del{( e_{i_1} \we \cdots \we e_{i_p})} = \sum_{k=1}^{p} 
 (-1)^{k-1}  e_{i_1} \we \cdots \we \widehat{e_{i_k}} \we 
        \cdots \we e_{i_p} $$ 
for $p \geq 2$ and $i_j \in [n]$.    
Let $M$ be a loopless matroid on $[n]$ with rank $\ell +1$.  

\begin{defn} 
The Orlik-Solomon algebra of $M$ 
is the quotient $A(M)$ of $E$ by the ideal $\ideal{\del M}$ 
generated by $\del{(e_{i_1} \we \cdots \we e_{i_s})}$ 
for every circuit $c=(i_1, \ldots, i_s)$ of $M$. 
\end{defn} 

If $1$ and $2$ are parallel, that is, $\{ 1, 2 \}$ is a circuit, 
then $e_1 = e_2$.  
So the Orlik-Solomon algebra of the simple matroid associated with $M$ 
is equal to that of $M$. 
The Orlik-Solomon algebra $A(M)$ has the natural grading. 
The linear map $\del$ on $E$ induces the linear map $\del_M$ on $A(M)$. 
Let $e_{\lambda} = \lambda_1 e_1 + \cdots + \lambda_n e_n \in E^1$. 
We call $\lambda = (\lambda_1, \ldots, \lambda_n)$ a weight of $M$. 
The left multiplication $e_{\lambda} \we : A^p(M) \to A^{p+1}(M)$ 
induces the complex $(A(M), e_{\lambda})$. 
Let $H(A(M), e_{\lambda})$ denote the cohomology of this complex. 
If $\lambda = 0$  then $H(A(M), e_{\lambda})$ is just $A(M)$, 
otherwise we have $H^0(A(M), e_{\lambda})=0$. 
If $\sum_{j=1}^{n} \lambda \not= 0$ then we have 
$H^p(A(M), e_{\lambda}) =0$ for all $p$ (see \cite{Yu1}). 
If $\del{e_{\lambda}}= \sum_{j=1}^{n} \lambda_j =0$ then 
$e_{\lambda}$ induces the complex $( \del_M{(A(M))}, e_{\lambda})$ 
and the cohomology $H( \del_M{(A(M))}, e_{\lambda})$, 
where  $\del_M{(A(M))}$ is the image of $\del_M$.  
It is known that 
$$ H^{p+1}(A(M), e_{\lambda})  = 
  H^{p+1}( \del_M{(A(M))}, e_{\lambda}) \oplus  
  H^p( \del_M{(A(M))}, e_{\lambda}).    $$ 
For a generic weight $\lambda$, 
Yuzvinsky \cite{Yu1} showed the vanishing theorem:    
$$   H^{k}( \del_M{(A(M))}, e_{\lambda}) = 0  
              \qquad \text{ for } k \not= \ell   $$ 
and hence we have 
$$ H^{k}(A(M), e_{\lambda})  = 0  
              \qquad \text{ for } k \not= \ell, \ell +1.  $$ 

An arrangement $\A$ of hyperplanes in $\mathbb{P}^{\ell}$ 
has the underlying matroid $M(\A)=M$ with rank $\ell +1$ as a combinatorial structure. 
The cohomology of the complement of $\A$ is isomorphic to 
 $\del_M{(A(M))}$  (see \cite{OT} and \cite{Ka}). 
If a weight $\lambda = (\lambda_i)_{i \in \A}$ satisfies 
some generic condition, 
then the cohomology of the complement of $\A$ 
with the coefficients in the rank one local system associated to $\lambda$ 
is isomorphic to  
$H( \del_M{(A(M))}, e_{\lambda})$ (see \cite{ESV, STV}). 
The local system cohomologys is an important subject in 
the multivariable theory of hypergeometric functions 
\cite{AK, OT2}.  
By the vanishing theorem \cite{Yu1}, 
for a generic weight $\lambda$, 
the local system cohomology vanishes 
in all but the top dimension.  
In this paper, our purpose is to construct matroids and arrangements 
with non-vanishing cohomology of Orlik-Solomon algebras, 
or rather $H^{\ell -1}(A(M), e_{\lambda}) \not= 0$.   

In particular, the case of $\ell =2$ was studied in \cite{Fa, LY}. 
Falk \cite{Fa} defined the \textit{resonance variety}, that is,  
the space of weights with non-vanishing cohomology. 
The resonance variety is deeply related to  
the cohomology support loci \cite{Ar} and  
the characteristic variety  \cite{Li, CS}. 
Libgober and Yuzvinsky \cite{LY} showed that, 
under some condition, weights with non-vanishing first cohomology 
are parametrized by Latin squares.  

In this paper, 
we prove that, in general, 
matroids associated to Latin hypercubes have 
weights with non-vanishing cohomology, 
by using decomposable relations arising from Latin hypercubes.  
This decomposable relation is the generalization of 
the relation discovered by Rybnikov (see \cite{Fa}).  
Moreover, in the case of $\ell =2$, 
we study well, using terms of Latin squares. 
In the last section, 
we shall give examples of realizations including the higher case. 
Some of them appear in the classical projective geometry 
(see Figure \ref{fig:pappus}, \ref{fig:kirkman} and \ref{fig:steiner}).  

We shall use the following notation and terminology.   
A $k$-set is a set with cardinality $k$. 
Denote the family of all $k$-subset of a set $S$ by $\binom{S}{k}$. 
Often, we regard a $p$-tuple $(i_1, \ldots, i_p)$ 
as a $p$-set $\{ i_1, \ldots, i_p \}$. 
We refer to \cite{Ox} for terminology of the matroid theory.

\section{ Non-vanishing Theorem } 

A Latin hypercube of dimension ${\ell}$ and order $m$ 
is an $m^{\ell}$-array  such that, 
if ${\ell}-1$ coordinates are fixed, 
the $m$ positions so determined contain a permutation of $m$ symbols. 
Let $K = [ k(i_1, \ldots, i_{\ell}) ]_{1 \leq i_1, \ldots, i_{\ell} \leq m}$ be 
a Latin ${\ell}$-dimensional hypercube on $[m]$, 
that is,  an $m^{\ell}$-matrix 
satisfying the condition 
$$  \{  k(i_1', i_2, \ldots, i_{\ell}) : i_1' \in [m]  \} = 
    \{  k(i_1, i_2', \ldots, i_{\ell}) : i_2' \in [m]  \} = \cdots  $$ 
$$ \cdots 
    = \{  k(i_1, i_2, \ldots, i_{\ell}') : i_{\ell}' \in [m]  \} = [m],  $$ 
for $1 \leq i_1, \ldots, i_{\ell} \leq m$. 
Define the family of $(\ell +1)$-subsets in $[n]$ associated to $K$ by   
$$  \C[K] = [ (i_1, m + i_2, 2m + i_3, \ldots, ({\ell}-1)m +i_{\ell}, 
  {\ell}m + k(i_1, \ldots, i_{\ell}) )  ]_{1 \leq i_1,\ldots, i_{\ell} \leq m}.    $$ 

On the other hand, 
a matroid is said to be \textit{$\ell$-generic} 
if it has no $i$-circuits for $i \leq \ell$. 
Note that 
an $1$-generic matroid is just a loopless matroid   
and  a $2$-generic matroid is just a simple matroid. 
The uniform matroid $U_{m,n}$ of rank $m$ is $m$-generic.  
So we can mention the main theorem as follows.  

\begin{thm}\label{thm:non-vanishing} 
Let $m \geq 2$,  $\ell \geq 2$ and $n=(\ell+1)m$. 
Let $K$ be a Latin $\ell$-dimensional hypercube on $[m]$.  
Then there exists a unique $\ell$-generic matroid $M[K]$ 
on $[n]$ with rank $\ell+1$,  
for which the family of all $(\ell+1)$-circuits is equal to $\C[K]$. 
This matroid has weights with non-vanishing cohomology, 
in fact, 
\begin{align*}
  H^{k}(A(M[K]), e_{\lambda}) &  =0 
          \qquad \text{ for }  k \leq \ell -2, \\ 
  H^{\ell -1}(A(M[K]), e_{\lambda}) & \not= 0, 
\end{align*}
for a non-zero weight 
$$  \lambda = 
( \underbrace{\lambda_1, \ldots, \lambda_1}_{m}, 
  \underbrace{\lambda_2, \ldots, \lambda_2}_{m}, 
    \cdots \cdots, 
  \underbrace{\lambda_{\ell +1}, \ldots, \lambda_{\ell +1}}_{m} );  
\quad \sum_{j=1}^{\ell +1} \lambda_j = 0.  $$ 
\end{thm} 

In the rest of this section, we will prove this theorem. 
First of all, we prove some lemmas.  

\begin{lem}\label{lem:r-circuit} 
A family $\C$ of $(\ell+1)$-subsets in $[n]$ satisfies the condition 
\begin{itemize} 
	\item[\cond{C_{\ell+1}}] 
if $C_1, C_2 \in \C$ and $\vert C_1 \cup C_2 \vert = \ell +2$ then 
every $(\ell+1)$-subset $C_3$ of $C_1 \cup C_2$ 
is a member of $\C$,  
\end{itemize} 
if and only if, 
there exists an $\ell$-generic matroid on $[n]$ for which 
the family of all $(\ell+1)$-circuits is equal to $\C$.  
\end{lem} 
\begin{proof} 
It is clear when $n< \ell+1$. 
Assume that $n \geq \ell +1$.
Let $\C$ be a family of $(\ell+1)$-subsets in $[n]$ 
satisfying \cond{C_{\ell+1}}. 
Let $I$ be a $\ell$-subset of $[n]$. 
Define $X_I = I \cup \{ e \in [n] : I \cup e \in \C \}$,  
$\binom{X_I}{\ell+1}   = \{ \text{all $(\ell+1)$-subsets of $X_I$.}  \}$, and 
$\binom{X_I}{\ell+1}_s  = \left\{ S \in \binom{X_I}{\ell+1} : 
                                  \card{S \setminus I} =s \right\}$. 
Note that $\binom{X_I}{\ell+1} = \cup_{s=1}^{\ell+1} \binom{X_I}{\ell+1}_s$. 
First of all, we show that $\binom{X_I}{\ell+1}_s$ is a subfamily of $\C$ 
by induction on $s$. 
For $s=1$, since $\binom{X_I}{\ell+1}_1 = \{ I \cup e \in \C \}$, it is clear. 
Let assume that $\binom{X_I}{\ell+1}_s \subset \C$ for $s \geq 1$.  
Take a member $S$ of $\binom{X_I}{\ell+1}_{s+1}$. 
Let $T := S \setminus I$ and $I' :=S \cap I$. 
Note that $S = I' \cup T$, 
$I' \subset I$, 
$T \subset X_I \setminus I$, 
$\card{I'}= \ell -s$ and 
$\card{T}= s+1$. 
Now we can choose $e \in I \setminus I'$ and 
$f_1, f_2 \in T$ with $f_1 \not= f_2$. 
By the inductive assumption, 
$C_1 := I' \cup e \cup ( T \setminus \{ f_1 \})$ 
and  $C_2 := I' \cup e \cup ( T \setminus \{ f_2 \})$ 
are in $\binom{X_I}{\ell+1}_{s} \subset \C$. 
We can check $C_1$ and $C_2$ satisfy the condition in \cond{C_{\ell+1}}, 
and $S$ is a $(\ell+1)$-subset of $C_1 \cup C_2$.  
So we have $S \in \C$. 
Therefore, we have $\binom{X_I}{\ell+1}_s \subset \C$ and 
hence $\binom{X_I}{\ell+1} \subset \C$.

Assume that $\C$ is not the family of all $(\ell+1)$-subsets of $[n]$. 
We shall show that 
$$ \I = \{ I \subset [n] : \vert I \vert \leq \ell+1, I \not\in \C \}  $$ 
is a matroid complex (see \cite{Ox}). 
Note that $I$ have all $i$-subsets of $[n]$ for $i < \ell+1$. 
Since 
$\emptyset \in \I$ and if $I' \subset I \in I$ then $I' \in \I$,   
we should prove 
the independence augmentation axiom for $\I$, 
that is, 
for $I_1, I_2 \in \I$ with $\card{I_2} = \card{I_1} +1$, 
there exists $e \in I_2 \setminus I_1$ such that $I_1 \cup \{ e \} \in \I$.  
If $\card{I_1} < \ell$, it is clear. 
Let $\card{I_1} = \ell$. 
Suppose that $I_1 \cup \{ e \} \not\in \I$  for all $e \in I_2 \setminus I_1$. 
Then we have $I_2 \subset X_{I_1}$.   
By the above claim,  we have $\binom{X_{I_1}}{\ell+1} \subset \C$ 
and hence we have $I_2 \in \C$, this is a contradiction.  
Therefore, $\I$ defines the matroid of rank $\ell+1$.  
The converse is easy by the circuit elimination axiom of the matroids 
(see \cite[1.1.4]{Ox}). 
\end{proof} 
\begin{rem} 
\begin{enumerate}
	\item 
When $\C = \emptyset$,  
the uniform matroid $U_{m,n}$  
of rank $m$ with $m \geq \ell+1$  is one of matroids in the above lemma. 
	\item 
If $\C$ consists of all $(\ell+1)$-subsets of $[n]$, 
the uniform matroid $U_{\ell,n}$ of rank $\ell$ 
is only one $\ell$-generic matroid in the above lemma. 
Otherwise, 
the rank of such a matroid is greater than $\ell$, and 
there exists uniquely such an $\ell$-generic matroid with rank $\ell +1$. 
\end{enumerate}
\end{rem} 

\begin{lem}\label{high-decom-rel}  
Let $n =({\ell}+1) m$. 
Let $a_s = e_{(s-1) m+1} + \cdots + e_{sm}$ for $1 \leq s \leq {\ell}+1$.  
For a Latin ${\ell}$-dimensional hypercube $K$ on $[m]$, 
we obtain the following decomposable relation 
\begin{align*} 
    \del{(a_1 \we a_2 \we \cdots \we a_{{\ell}+1})}      
 & = - (a_1 - a_2) \we (\del{(a_2 \we \cdots \we a_{{\ell}+1})})  \\ 
 & = (-1)^{{\ell}} m  (a_1 - a_2) \we (a_2 - a_3) \we 
            \cdots \we (a_{\ell} - a_{{\ell}+1})  \\ 
 & = m \sum_{S \in \C[K]}   \del{(e_S)}, 
\end{align*}  
where $e_S = e_{i_1} \we \cdots \we e_{i_p}$ 
for a $p$-tuple $(i_1, \ldots, i_p)$. 
\end{lem} 
\begin{proof} 
The first and second equations are obtained by 
\begin{align*} 
   &  \del{(a_1 \we a_2 \we \cdots \we a_{{\ell}+1})}      
 = \del{( (a_1- a_2) \we a_2 \we \cdots \we a_{{\ell}+1})}    \\   
 & =  \del{(a_1- a_2)} \we a_2 \we \cdots \we a_{{\ell}+1}  
       - (a_1 - a_2) \we (\del{(a_2 \we \cdots \we a_{{\ell}+1})})   \\ 
 & =   - (a_1 - a_2) \we (\del{(a_2 \we \cdots \we a_{{\ell}+1})}) = \cdots   \\ 
 & =  (-1)^{\ell}   (a_1 - a_2) \we (a_2 - a_3) \we \cdots \we (a_{\ell} - a_{{\ell}+1}) 
     \we \del{(a_{{\ell}+1})}   \\ 
 & = (-1)^{\ell} m (a_1 - a_2) \we (a_2 - a_3) \we \cdots \we (a_{\ell} - a_{{\ell}+1}). 
\end{align*}  

Let $E_s = \{ (s-1) m + 1, (s-1) m + 2, \ldots, sm  \}$ for $1 \leq s \leq {\ell}+1$. 
Note that $E_1 \cup \cdots \cup E_{{\ell}+1} = [n]$. 
We regard $K$ as a Latin hypercube $\tilde{K} = (\tilde{k}(i_1, \ldots, i_{\ell}) )$ 
with $s$-axis indexed by $E_s$ for $1 \leq s \leq {\ell}$ 
and symbol set $E_{{\ell}+1}$.  
We note that 
$\tilde{k}(i_1, \ldots, i_{\ell}) = {\ell}m + k(i_1, \ldots, i_{\ell}) 
\in E_{{\ell}+1}$. 
Since  
$\del{(e_{1} \we  \cdots \we e_{k} \we e_{k+1})} 
      = \del{(e_{1} \we  \cdots \we e_k)} \we e_{k+1} 
       + (-1)^{k} e_{1} \we \cdots \we e_k$,  
we have 
$$ (-1)^{k} e_{1} \we e_{2} \we \cdots \we e_k  
       = - \del{(e_{1} \we  \cdots \we e_k)} \we e_{k+1} 
         +  \del{(e_{1} \we  \cdots \we e_k \we e_{k+1})}.   $$ 
Hence, we can get   
\begin{align*} 
  & (-1)^{\ell} m \cdot a_1 \we \cdots \we a_{\ell} = 
  m \sum_{i_1 \in E_1, \ldots, i_{\ell} \in E_{\ell}} (-1)^{\ell}  
     e_{i_1} \we  \cdots \we e_{i_{\ell}}          
   = m \times \\ 
 & \sum_{i_1 \in E_1, \ldots, i_{\ell} \in E_{\ell}} \left\{ 
 - \del{(e_{i_1} \we  \cdots \we e_{i_{\ell}})} \we e_{\tilde{k}(i_1, \ldots, i_{\ell})} 
 + \del{(e_{i_1} \we \cdots \we e_{i_{\ell}} \we e_{\tilde{k}(i_1, \ldots, i_{\ell})})} 
   \right\}.  
\end{align*} 
The second term is 
$$  \sum_{i_1 \in E_1, \ldots, i_{\ell} \in E_{\ell}}  
  \del{(e_{i_1} \we \cdots \we e_{i_{\ell}} \we e_{\tilde{k}(i_1, \ldots, i_{\ell})})}   
  = \sum_{S \in \C[K]}   \del{(e_S)}. $$ 
On the other hand, since $K$ is  a Latin hypercube, 
we have  
\begin{align*} 
 & \sum_{i_1 \in E_1, \ldots, i_{\ell} \in E_{\ell}}  
 \del{(e_{i_1} \we \cdots \we e_{i_{\ell}})} \we e_{\tilde{k}(i_1, \ldots, i_{\ell})} \\ 
 & = \sum_{i_1, \ldots, i_{\ell}}  
   \left( \sum_{p =1}^{\ell} (-1)^{p-1}  
   e_{i_1} \we \cdots \we \widehat{e_{i_p}} \we  \cdots \we e_{i_{\ell}} \right) 
    \we  e_{\tilde{k}(i_1, \ldots, i_{\ell})}    \\ 
 & = \sum_{p =1}^{{\ell}}  \sum_{i_1, \ldots, i_{\ell}}  
   \left( (-1)^{p-1}  
   e_{i_1} \we \cdots \we \widehat{e_{i_p}} \we  \cdots \we e_{i_{\ell}} \right) 
    \we  e_{\tilde{k}(i_1, \ldots, i_{\ell})}    \\ 
 & = \sum_{p =1}^{\ell}  \sum_{i_1, \ldots, \widehat{i_p}, \ldots, i_{\ell}}  
   \left( (-1)^{p-1}  
   e_{i_1} \we \cdots \we \widehat{e_{i_p}} \we  \cdots \we e_{i_{\ell}} \right) 
    \we  \sum_{i_p} e_{\tilde{k}(i_1, \ldots, i_{\ell})}    \\ 
 & = \sum_{p =1}^{\ell}  \sum_{i_1, \ldots, \widehat{i_p}, \ldots, i_{\ell}}  
   \left( (-1)^{p-1}  
   e_{i_1} \we \cdots \we \widehat{e_{i_p}} \we  \cdots \we e_{i_{\ell}} \right) 
    \we  a_{{\ell}+1},      
\end{align*} 
and 
\begin{align*} 
 \del{(a_1 \we \cdots \we a_{\ell})}  
    & = \del{(a_p)} \sum_{p =1}^{\ell} (-1)^{p-1}  
   a_1 \we \cdots \we \widehat{a_p} \we  \cdots \we a_{\ell}  \\   
  & = m \sum_{p =1}^{\ell} (-1)^{p-1}  
       \sum_{i_1, \ldots, \widehat{i_p}, \ldots, i_{\ell}}  
    e_{i_1} \we \cdots \we \widehat{e_{i_p}} \we  \cdots \we e_{i_{\ell}}.   
\end{align*} 
Therefore we obtain  
$$  (-1)^{\ell} m \cdot a_1 \we \cdots \we a_{\ell} = 
   - \del{(a_1 \we \cdots \we a_{\ell})} \we  a_{{\ell}+1}  
   + m  \sum_{S \in \C[K]}   \del{(e_S)}  $$ 
and hence we have 
\begin{align*} 
 \del{(a_1 \we \cdots \we a_{\ell} \we a_{{\ell}+1})}  & = 
   \del{(a_1 \we \cdots \we a_{\ell})} \we  a_{{\ell}+1}  + 
      (-1)^{\ell} m \cdot a_1 \we \cdots \we a_{\ell} \\ 
      & =   m  \sum_{S \in \C[K]}   \del{(e_S)}. 
\end{align*} 
\end{proof} 

\begin{proof}[Proof of Theorem \ref{thm:non-vanishing}]  
Let $K$ be a Latin $\ell$-dimensional hypercube on $[m]$.  
By the construction of $\C[K]$, 
for $C_1, C_2 \in \C[K]$ with $C_1 \not= C_2$, 
we have $\vert C_1 \cap C_2 \vert = \ell -1$ and 
$\vert C_1 \cup C_2 \vert = \ell +3$. 
Hence, due to Lemma \ref{lem:r-circuit} and its remark,  
there exists 
a unique  $\ell$-generic matroid $M[K]$ with rank $\ell +1$. 
In general, 
for an $\ell$-generic matroid $M$ and a non-zero weight $\lambda$ of $M$, 
we have $H^k(M, e_{\lambda})=0$ for $k \leq \ell -2$. 
Thus, we shall prove 
$H^{\ell -1}(A(M[K]), e_{\lambda}) \not= 0$.  
Let $\lambda$ be a weight given in the statement, 
and assume without loss of generality that $\lambda_1 \not= 0$. 
Since $\sum_{j=1}^{\ell +1} \lambda_j = 0$, we have 
\begin{align*}
  e_{\lambda}  & = 
   \lambda_1 (e_1 + \cdots + e_m )  + \cdots + 
    \lambda_{\ell +1} (e_{\ell m+1} + \cdots + e_{(\ell +1)m} )   \\ 
  & =    \lambda_1 a_1  +  \lambda_2 a_2  + \cdots + 
    \lambda_{\ell +1} a_{\ell +1}   \\ 
  & =  \lambda_1 (a_1 - a_2) 
    + (\lambda_1 + \lambda_2) (a_2 - a_3)  + \cdots  
    + (\lambda_1 + \cdots + \lambda_{\ell})  (a_{\ell} - a_{\ell +1}), 
\end{align*}
where $a_j$ is defined in Lemma \ref{high-decom-rel}. 
Define a $(\ell -1)$-form 
\begin{align*}
 b & := \del{(a_2 \we a_3 \we \cdots \we a_{{\ell}+1})} \\ 
 & = (-1)^{{\ell -1}} m  (a_2 - a_3) \we (a_3 - a_4) \we 
            \cdots \we (a_{\ell} - a_{{\ell}+1}). 
\end{align*}
By Lemma \ref{high-decom-rel}, we have  
$$   e_{\lambda} \we b = \lambda_1 (a_1 - a_2) \we  b  
  \in \ideal{\del{M[K]}}, $$ 
that is, $e_{\lambda} \we b$ vanishes 
in the Orlik-Solomon algebra $A(M[K])$.  
Since $M[K]$ is $\ell$-generic,  
the $(\ell -1)$-form $b$ is not in $\ideal{\del{M[K]}}$. 
Finally, we shall check that 
$b$ is a non-vanishing cohomology class 
in $H^{\ell -1}(M[K], e_{\lambda})$.  

For a finite set $\{ e_1, \ldots, e_n \}$, 
denote $E(e_1, \ldots, e_n)$ 
the graded exterior algebra over $R$ generated 
by $1$ and degree-one elements $e_1, \ldots, e_n$.  
Note that $E(e_2, \ldots, e_n)$ is 
a subalgebra of $E(e_1, \ldots, e_n)$. 
Let $e_{\lambda} = \lambda_1 e_1 + \cdots +  \lambda_n e_n$ 
with $\lambda_i \in R$ and $\lambda_1 \not= 0$.  
Then we have 
$E(e_1, \ldots, e_n) = E(e_{\lambda}, e_2, \ldots, e_n)$.   
It is easy to see the following: 
if $\omega \in E(e_2, \ldots, e_n)$ with $\omega \not= 0$, 
then $\omega$ is not belong to the ideal of $E(e_1, \ldots, e_n)$ 
generated by $e_{\lambda}$. 

By the above,  
since $b$ is in $E(e_{m+1}, \ldots, e_{n})$ and $\lambda_1 \not= 0$, 
$b$ is not in the ideal of $E(e_1, \ldots, e_n)$ 
generated by $e_{\lambda}$, that is, 
there exists no $(\ell -2)$-form $\eta$ 
with $e_{\lambda} \we \eta = b$. 
This completes the proof. 
\end{proof}

\section{ The case of $\ell=2$  } 

We refer to \cite{CH} for the Latin squares.   
A \textit{Latin square} of order $m$  
is a Latin hypercube of dimension 2 and order $m$, 
that is,  an $m \times m$ matrix 
with entries in an $m$-set 
(we call the \textit{symbol set}.)  
such that each element occurs 
exactly once in each row  and 
exactly once in each column.   
The two Latin squares $K$ and $K'$ are 
\textit{isotopic} if 
$K'$ is obtained by 
permutations of rows, permutations of columns, 
and a bijection from the symbol set of $K$.  
Let $E_1$, $E_2$ and $E_3$ be three $m$-sets and 
let $K$ be a Latin square 
with rows indexed by $E_1$, columns by $E_2$, and symbols by $E_3$. 
Define  
$T(K) = \{ \{ x_1, x_2, x_3 \} : x_i \in E_i (i=1,2,3), 
k_{x_1,x_2} = x_3 \}$. 
For any permutation $\sigma$ of $\{ 1,2,3 \}$, 
the $\sigma$-\textit{conjugate} of $L$ 
is the Latin square $K_{\sigma}$ 
with rows indexed by $E_{\sigma 1}$,  columns by $E_{\sigma 2}$, 
and symbols by $E_{\sigma 3}$,  
defined by $T(K) = T(K_{\sigma})$.  
The two Latin squares $K$ and $K'$ are 
\textit{main class isotopic} if 
$K'$ is isotopic to any conjugate of $K$. 

Let $K = ( k_{i,j} )$ be a a Latin square on $[m]$, 
that is, an $m \times m$-matrix 
satisfying the condition 
$\{ k_{i,1}, k_{i,2}, \ldots, k_{i,m} \} 
   = \{ k_{1,j}, k_{2,j}, \ldots, k_{m,j} \} = [m]$    
for $1 \leq i, j \leq m$.    
In the previous section, 
we define $\C{[K]}$ by the family 
$$  
\begin{bmatrix}
  (1, m+1, 2m+k_{1,1}) &  (1, m+2, 2m+k_{1,2}) & \cdots & (1, 2m, 2m+k_{1,m}) \\ 
  (2, m+1, 2m+k_{2,1}) &  (2, m+2, 2m+k_{2,2}) & \cdots & (1, 2m, 2m+k_{2,m}) \\ 
             \vdots &                    &        &    \vdots   \\ 
  (m, m+1, 2m+k_{m,1}) &  (m, m+2, 2m+k_{m,2}) & \cdots & (1, 2m, 2m+k_{m,m}) 
\end{bmatrix}.   $$ 
We can see $K$  
as a Latin square $\tilde{K}$ 
with rows indexed by $\{ 1,2,\ldots,m \}$,  
columns by $\{ m+1,m+2,\ldots,2m \}$, 
and symbols by $\{ 2m+1, 2m+2, \ldots, 3m \}$.   
So we can consider $\C[K]=T(\tilde{K})$.  
By Theorem \ref{thm:non-vanishing}, 
there exists a unique simple matroid $M[K]$ 
on $[n]$ with rank $3$,  
for which the family of all $3$-circuits is equal to $\C[K]$. 
The simple matroid $M[K]$ has 
weights with non-vanishing first cohomology. 
\begin{prop} 
Let $m \geq 2$.   
If $K_1$ and $K_2$ are main class isotopic Latin squares 
then matroids $M[K_1]$ and $M[K_2]$ are isomorphic. 
If a Latin square $K_1$ is not main class isotopic to $K_2$ 
then matroids $M[K_1]$  is not isomorphic to $M[K_2]$. 
\end{prop} 
\begin{proof} 
It is clear by the definition of main class isotopic Latin squares. 
\end{proof} 
\begin{rem} 
The number of main class isotopic Latin squares of order $m \leq 8$ is known 
(see \cite{CH}).   
\begin{center}
\begin{tabular}{|c||c|c|c|c|c|c|c|c|} 
\hline   $m =$ & $1$ & $2$ & $3$ & $4$ & $5$ & $6$ & $7$ & $8$ \\ 
\hline  main classes & $1$ & $1$ & $1$ & $2$ & $2$ & $12$ & $147$ & $283,657$  \\  
\hline 
\end{tabular} 
\end{center}
\end{rem} 

The two Latin squares $K=(k_{i,j})$ and $K'=(k'_{i,j})$ 
of same order are 
\textit{orthogonal} if 
all pairs $(k_{i,j}, k'_{i,j})$ are distinct. 
A set of Latin squares of order $m$ is 
\textit{mutually orthogonal} 
if any two distinct squares are orthogonal. 

\begin{thm}\label{thm:orth}  
Let $m \geq 1$, $s \geq 1$ and $n=(s+2)m$.  
Let $K_1$, \ldots, $K_s$ be mutually orthogonal Latin squares on $[m]$.  
Then there exists 
a simple matroid $M[K_1, \ldots, K_s]$ on $[n]$ satisfying 
$$ \dim{H^1(A(M[K_1, \ldots, K_s]), e_{\lambda})} = s  $$ 
for a non-zero weight 
$$  \lambda = 
( \underbrace{\lambda_1, \ldots, \lambda_1}_{m}, 
  \underbrace{\lambda_2, \ldots, \lambda_2}_{m}, 
    \cdots \cdots, 
  \underbrace{\lambda_{s+2}, \ldots, \lambda_{s+2}}_{m} );  
\quad \sum_{j=1}^{s+2} \lambda_j = 0.  $$ 
\end{thm} 

\begin{proof} 
By Lemma \ref{lem:r-circuit} in the case of $\ell =2$,   
a family $\C$ of $3$-subsets in $[n]$ satisfies the condition 
\begin{itemize} 
	\item[\cond{C_3}] 
if $\{ i,j,k \}$ and $\{ i,j,l \}$ are members of $\C$ 
then $\{ i,k,l \}$ and $\{ j,k,l \}$ are members of $\C$,  
\end{itemize} 
if and only if, 
there exists a simple matroid on $[n]$ for which 
the family of all $3$-circuits is equal to $\C$.  
Recall that 
the set of flats of a matroid is a geometric lattice.  
The closure of $C \in \C$ is the set 
$\cup \{ C' \in \C : \vert C' \cap C \vert \geq 2 \}$, 
that is a flat of rank 2. 
A $2$-subset contained in no $C \in \C$ 
is a flat of rank 2.  

\textbf{Construction of $M[K_1, \ldots, K_s]$}: 
Let $K_1$, \ldots, $K_s$ be mutually orthogonal Latin squares on $[m]$.  
A sift Latin square 
$\tilde{K_p} =(\tilde{k}_{i,j}^p)$ 
associated to $K_p=(k_{i,j}^p)$ 
is given by a Latin square 
with row indexed by $\{ 1,2,\ldots,m \}$,  
column by $\{ m+1,m+2,\ldots,2m \}$, 
and symbols by $\{ (p+1)m+1,(p+1)m+2,\ldots, (p+2)m \}$,  
given by $\tilde{k}_{i,j}^p = (p+1)m + k_{i,j}^p$ 
for $1 \leq i \leq m$ and $m+1 \leq j \leq 2m$.  
We define  
\begin{align*}
  & \C[K_1, \ldots, K_s] : = 
   T(\tilde{K_1}) \cup \cdots \cup T(\tilde{K_s}),    \\ 
  & X_{i,j}  := \{ i, j, \tilde{k}_{i,j}^1, \ldots, \tilde{k}_{i,j}^s \} 
     \quad \text{ for } 1 \leq i \leq m, m+1 \leq j \leq 2m,  \  \text{ and }  \\ 
  & \C  := \C[K_1, \ldots, K_s] \cup 
 \left( \bigcup_{1 \leq i \leq m, m+1 \leq j \leq 2m} \binom{X_{i,j}}{3} \right).     
\end{align*} 
By mutually orthogonality, 
we have $\vert C \cap X_{i,j} \vert =1$ 
for any $C \in \C[K_1, \ldots, K_s]$ not contained in $X_{i,j}$,  
and 
$\vert X_{i,j} \cap X_{k,l} \vert =1$ 
for $(i,j) \not= (k,l)$. 
This implies that $\C$ satisfies \cond{C_3}.  
If $m \geq 2$ 
then we obtain a simple matroid $M[K_1, \ldots, K_s]$ on $[n]$ with rank $3$ 
such that $\C$ is the family of all $3$-circuits. 
If $m=1$ then $\C$ gives the uniform matroid $U_{2,n}$. 

\textbf{Non-vanishing}: 
Let $a_i = e_{(i-1)m} + e_{(i-1)m +1} + \cdots + e_{(i-1)m}$ 
for $i=1,2, \ldots, s+2$. 
By Lemma \ref{high-decom-rel}, we have 
$$ (a_1 - a_i) \we (a_2 - a_i) \in \ideal{\del{M[K_1, \ldots, K_s]}} $$ 
for $3 \leq i \leq s+2$. 
We take two one-forms 
$$ e_{\lambda^t} = 
   \lambda_1^t a_1 + \lambda_2^t a_2 + \cdots 
 + \lambda_{s+2}^t a_{s+2}   $$ 
with $\sum_{j=1}^{s+2} \lambda_j^t = 0$ for $t=1,2$. 
Since 
$e_{\lambda^1} = 
   \lambda_2^1 (a_2 -a_1) + \cdots 
 + \lambda_{s+2}^1 (a_{s+2}-a_1)$ and  
$e_{\lambda^2} = 
   \lambda_1^2 (a_1 -a_2) +  \cdots 
 + \lambda_{s+2}^2 (a_{s+2}-a_2)$, 
we have 
$e_{\lambda^1} \we e_{\lambda^2} 
\in \ideal{\del{M[K_1, \ldots, K_s]}}$. 
This implies 
$\dim{H^1(A(M[K_1, \ldots, K_s]), e_{\lambda})} = s$. 
\end{proof} 
\begin{rem} 
When $m=1$,  the matroid in this theorem is the uniform matroid $U_{2,n}$ with rank $2$. 
When $m \geq 2$, the matroid $M[K_1, \ldots, K_s]$ has rank $3$. 
\end{rem} 
\begin{rem} 
There exists a Latin square of order $m$ for $m \geq 1$. 
Let $N(m)$ be the maximum number of 
mutually orthogonal Latin squares of order $m$. 
The following is known (see \cite{CH}).  
\begin{itemize}
	\item  $N(0)=N(1) = \infty$ and  $1 \leq N(m) \leq m-1$ for every $m>1$. 
	\item  
If $m$ is a prime power then $N(m)=m-1$. 
	\item  
If $m \not\equiv 2 \mod{4}$, then $N(m) \geq 2$. 
	\item  $N(p \times q) \geq \min{\{ N(p), N(q) \}}$.  
	\item  $N(2)=1$, $N(3)=2$, $N(4)=3$, $N(5)=4$, 
           $N(6)=1$, $N(7)=6$, $N(8)=7$. 
\end{itemize}
\end{rem} 
\begin{rem}\label{rem:3}  
In the case of $s=1$,  we have 
$\dim{H^1(A(M[K]), e_{\lambda})} = 1$ 
for non-zero one-form 
$$ e_{\lambda} = 
   \lambda_1 (e_1 + \cdots + e_m ) 
 + \lambda_2 (e_{m+1} + \cdots + e_{2m} ) 
 + \lambda_{3} (e_{2m+1} + \cdots + e_{3m} ) $$ 
with $\lambda_1 + \lambda_2 + \lambda_3 = 0$. 
\end{rem}  

Let $M$ and $M'$ be loopless matroids $M$ on $[n]$ of rank 3. 
We call $M'$ a \textit{degeneration} of $M$ 
if the family of $3$-circuits of $M'$ contains that of $M$.    
Mostly, 
degenerations of $M[K_1, \ldots, K_s]$ have weights with non-vanishing first cohomology. 
The uniform matroid $U_{2,n}$ of rank $2$ is its degeneration.    
Next, without $U_{2,n}$,  we shall construct its degeneration 
with non-vanishing first cohomology. 

\begin{prop}\label{prop:deg1}  
Let $m \geq 2$, $s \geq 1$ and $n=(s+2)m$.  
Let $K_1$, \ldots, $K_s$ be mutually orthogonal Latin squares on $[m]$.  
Let $M_i$ be a simple matroid on 
$I_i := \{ (i-1)m+1, (i-1)m+2, \ldots, im  \}$ for $i=1,2, \ldots, s+2$. 
There exists  
a simple matroid $M[K_1, \ldots, K_s  : M_1, \ldots ,M_{s+2} ]$ 
with rank $3$ 
such that it is a degeneration of $M[K_1, \ldots, K_s]$ 
and its restriction on $I_i$ is $M_i$ for $i=1,2, \ldots, s+2$. 
Then we have 
$$ \dim{H^1(A(M[K_1, \ldots, K_s  : M_1, \ldots ,M_{s+2}]), 
    e_{\lambda})} = s $$ 
for a weight ${\lambda}$ given in Theorem \ref{thm:orth}. 
\end{prop} 
\begin{proof} 
Let $\C_3{(M_1, \ldots ,M_{s+2})}$ be 
the union of families of $3$-curcuits of $M_i; i=1, \ldots, s+2$. 
For a $3$-curcuit $C_i$ of $M_i$  
and $C \in \C[K_1, \ldots, K_s]$,  
we have $C_i \cap C_j = \emptyset$ for $i \not= j$ 
and $\vert C_i \cap C \vert =1$. 
Thus  $\C[K_1, \ldots, K_s] \cup \C_3{(M_1, \ldots ,M_{s+2})}$ satisfies \cond{C_3} and 
it yields  a simple matroid 
$M[K_1, \ldots, K_s  : M_1, \ldots ,M_{s+2} ]$ in this statement. 
By the same argument as that in the proof of Theorem \ref{thm:orth}, 
we can prove the proposition. 
\end{proof} 

\begin{rem} 
A realization of $M[K_1, \ldots, K_s  : M_1, \ldots ,M_{s+2} ]$  
is a $(s+2,m)$-net in $\mathbb{P}^2$ defined in \cite{Yu3}.  
Therefore, 
there is no $(k,m)$-net for $k > N(m) +2$. 
In particular, there is no $(k,6)$-net for $k > 3$. 
\end{rem} 

In a Latin square $K$,  
a $s \times s$-matrix obtained by $s$ rows and $s$ columns 
is called a Latin $s$-\textit{subsquare} of $K$  
if it forms a Latin square of order $s$. 
Let $K$ be a Latin square on $[m]$ and $J$ be a subsquare of $K$. 
We treat $\tilde{J}$ as a subsquare of $\tilde{K}$. 
$\tilde{J}$ has row index set $I_1(J)$, 
column index set $I_2(J)$ and symbol set $I_3(J)$ 
where $I_1(J) \subset I_1$, 
$I_2(J) \subset I_2$, 
$I_3(J) \subset I_3$ 
and $\vert I_1(J) \vert = \vert I_2(J) \vert = \vert I_3(J) \vert$. 
We define 
$X{(J)} = I_1(J) \cup I_2(J) \cup I_3(J)$.  
\begin{prop}\label{prop:deg2} %
Let $J$ be a subsquare of a Latin square $K$ on $[m]$.    
There exists a simple degeneration $M[K; J]$ of $M[K]$,  
whose restriction on $X{(J)}$ is the uniform matroid of rank 2.  
Then we have 
$$  \dim{H^1(A(M[K; J]), e_{\lambda})} = 1  $$  
for a weight $\lambda$ given in Remark \ref{rem:3}. 
\end{prop} 

\begin{proof} 
Let $\C = \C[K] \cup \binom{X{(J)}}{3}$. 
Since $J$ is a subsquare of $K$, 
for $C \in \C[K] \setminus  \binom{X{(J)}}{3}$, 
we have $\vert C \cap X{(J)} \vert =1$. 
This leads to \cond{C_3} for $\C$. 
By the same way of Proposition \ref{prop:deg1}, we can show this. 
\end{proof} 
\begin{rem} 
The following is known (see \cite{CH}). 
\begin{itemize}
	\item 
There exists a Latin square of order $m$ with a proper $k$-subsquare 
if and only if $k \leq \left[ \frac{m}{2} \right]$. 
	\item 
There exists a Latin square of order $m$ with no proper subsquares 
if $m \not= 2^a3^b$ or if $m= 3$,$9$,$12$,$16$,$18$,$27$,$81$ or $243$. 
\end{itemize}
\end{rem} 

We note that  
there are degenerations of matroids associated to Latin square 
with non-vanishing cohomology,   
except for those of Proposition \ref{prop:deg1} and \ref{prop:deg2}. 
Especially, it is not necessary to be simple, 
for example, see Section \ref{arr:deg}.  

\section{ Arrangements } 

For a matroid $M$, 
an arrangement over a field $F$ with underlying matroid $M$ 
is called a $F$-\textit{realization} or \textit{representation} of $M$. 
A matroid is said to be \textit{realizable} or \textit{representable} 
over $F$ if $M$ has a $F$-realization.  
We shall find realizations of matroids obtained in the previous section.  
In this section, we will know the following: 
\begin{prop} 
If $1 \leq m \leq 4$ then the matroid $M[K]$ 
associated to a Latin square $K$ on $[m]$ 
is realizable over real. 
\end{prop} 
In addition, 
these realizations are arrangements 
appearing in the classical projective geometry 
(Figure \ref{fig:pappus}, \ref{fig:kirkman} and \ref{fig:steiner}).  
Besides, we shall give many other examples including the higher case. 

\subsection{$m=1$} 
Lemma \ref{high-decom-rel} implies 
$(e_1 - e_3) \we (e_2 - e_3) = \del{(e_1 \we e_2 \we e_3)}$. 
The matroid $M[K]$ is realized by 
the arrangement in $\mathbb{P}^2$ 
consisting of three lines through one point. 

\subsection{$m=2$ (Falk \cite{Fa})} 
We have only one main class isotopic Latin square 
$K = 
\begin{pmatrix} 
 1 & 2 \\ 
 2 & 1 
\end{pmatrix}$. 
The decomposable relation is 
$(e_1 + e _2 - e_5 - e_6) \we (e_3 + e_4 - e_5 -e_6) 
= \del{(e_1 \we e_3 \we e_5)} +  \del{(e_1 \we e_4 \we e_6)} 
   +  \del{(e_2 \we e_3 \we e_6)} +  \del{(e_2 \we e_4 \we e_5)}$. 
The matroid $M[K]$ is realized by 
the arrangement in $\mathbb{P}^2$  arising from the Ceva Theorem 
(the left side in Figure \ref{fig:pappus}). 

\subsection{$m=3$} 
\begin{figure}[h] 
	\includegraphics[width=.3\textwidth]{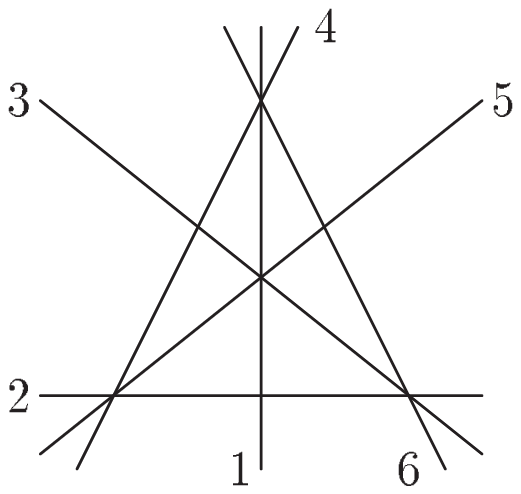} %
\qquad \quad  
	\includegraphics[width=.4\textwidth]{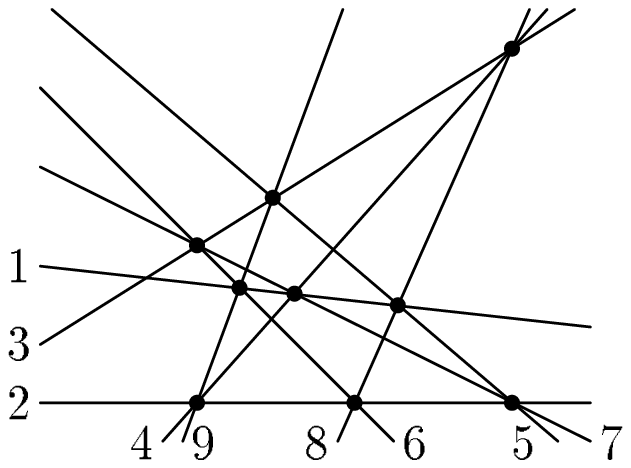} %
\caption{The Ceva Theorem and the Pappus Theorem}
\label{fig:pappus}
\end{figure}

We have only one main class isotopic Latin square, 
which is given by  
$$  
K = 
\begin{pmatrix} 
 1 & 2 & 3 \\ 
 3 & 1 & 2 \\ 
 2 & 3 & 1
\end{pmatrix}, \  
\C[K] = 
\begin{bmatrix}
  (1,4,7) & (1,5,8) & (1,6,9) \\ 
  (2,4,9) & (2,5,7) & (2,6,8) \\ 
  (3,4,8) & (3,5,9) & (3,6,7)  
\end{bmatrix}.  $$ 
The realization is given by 
the arrangement of 9 lines in $\mathbb{P}^2$ 
arising from the Pappus Theorem (the right side in Figure \ref{fig:pappus}).  

\subsection{$m=4$}\label{arr:m=4} 
\begin{figure}[h] 
	\includegraphics[width=.7\textwidth]{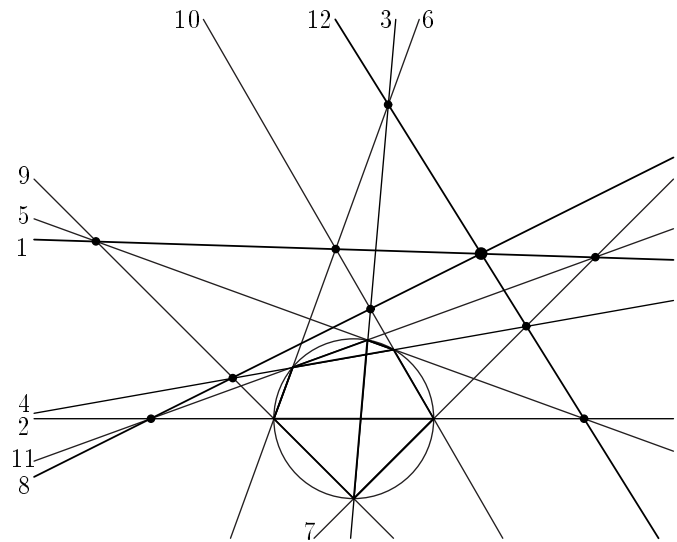} %
\caption{The Kirkman Theorem}
\label{fig:kirkman}
\end{figure}

\begin{figure}[h] 
	\includegraphics[width=.7\textwidth]{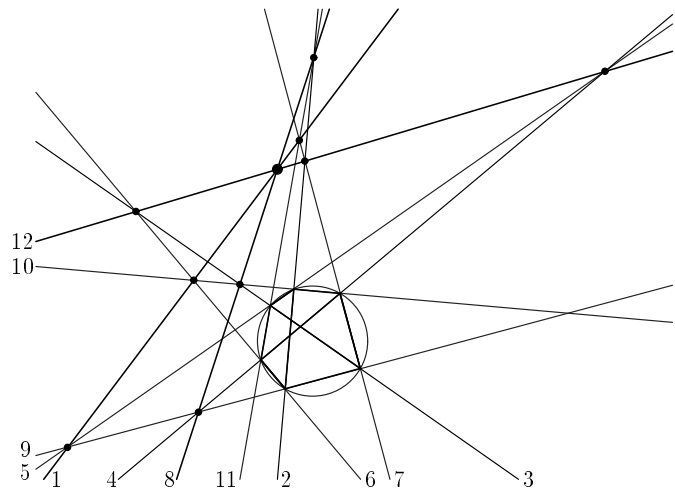} %
\caption{The Steiner Theorem}
\label{fig:steiner}
\end{figure}

There are two main class isotopic Latin squares, 
that we can give by 
$$  
K_1 = 
\begin{pmatrix} 
 1 & 2 & 3 & 4 \\ 
 4 & 1 & 2 & 3 \\ 
 3 & 4 & 1 & 2 \\ 
 2 & 3 & 4 & 1 \\  
\end{pmatrix}, \  
\C[K_1] = 
\begin{bmatrix}
  (1,5,9)  & (1,6,10) & (1,7,11) & (1,8,12)  \\ 
  (2,5,12) & (2,6,9)  & (2,7,10) & (2,8,11)  \\ 
  (3,5,11) & (3,6,12) & (3,7,9)  & (3,8,10)  \\ 
  (4,5,10) & (4,6,11) & (4,7,12) & (4,8,9)   
\end{bmatrix}, $$ 
$$  
K_2 = 
\begin{pmatrix} 
 1 & 2 & 3 & 4 \\ 
 2 & 1 & 4 & 3 \\ 
 3 & 4 & 1 & 2 \\ 
 4 & 3 & 2 & 1 \\  
\end{pmatrix}, \ 
\C[K_2] = 
\begin{bmatrix}
  (1,5,9)  & (1,6,10) & (1,7,11) & (1,8,12)  \\ 
  (2,5,10) & (2,6,9)  & (2,7,12) & (2,8,11)  \\ 
  (3,5,11) & (3,6,12) & (3,7,9)  & (3,8,10)  \\ 
  (4,5,12) & (4,6,11) & (4,7,10) & (4,8,9)   
\end{bmatrix}. $$ 
The matroid $M[K_1]$ or $M[K_2]$ 
is realized by the arrangement of 12 lines in $\mathbb{P}^2$ 
defined by Figure \ref{fig:kirkman} or \ref{fig:steiner}, 
which is arising from 
the Kirkman Theorem or the Steiner Theorem, respectively 
 (see \cite[Chapter 16]{Pr}).  

\subsection{Degenerations}\label{arr:deg} 
Let $K_1$ and $K_2$ be in the preceding section.  

\begin{figure}[h] 
	\includegraphics[width=.4\textwidth]{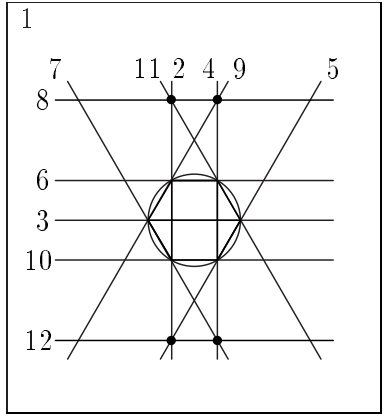} %
\caption{Degeneration of Kirkman's arrangement}
\label{fig:dkirkman}
\end{figure}

Let $J$ be a subsquare of $K_1$ given by 
$$  J = 
\begin{pmatrix} 
   & 2 &   & 4 \\ 
   &   &   &   \\ 
   & 4 &   & 2 \\ 
   &   &   &   \\  
\end{pmatrix}.    $$  
By Proposition \ref{prop:deg2}, 
we obtain $X{(J)} = \{ 1,3,6,8,10,12 \}$ and the matroid $M[K_1;J]$.   
Let $M_1$ be a simple matroid on $[4]$ for which 
the family of $3$-circuits is $\{ (1,2,4) \}$. 
By Proposition \ref{prop:deg1}, we have the matroid $M[K_1;M_1]$. 
Furthermore, the family 
$\C[K_1] \cup \binom{X(J)}{3} \cup \C_3(M_1)$ satisfies \cond{C_3} 
and then yields the matroid $M[K_1:M_1;J]$ 
with non-vanishing first cohomology.   
This matroid $M[K_1:M_1;J]$ is realized by 
the arrangement of 
11 lines in $\mathbb{C}^2$ with the infinite line $1$ 
in Figure \ref{fig:dkirkman}.   

\begin{figure}[h] 
	\includegraphics[width=.4\textwidth]{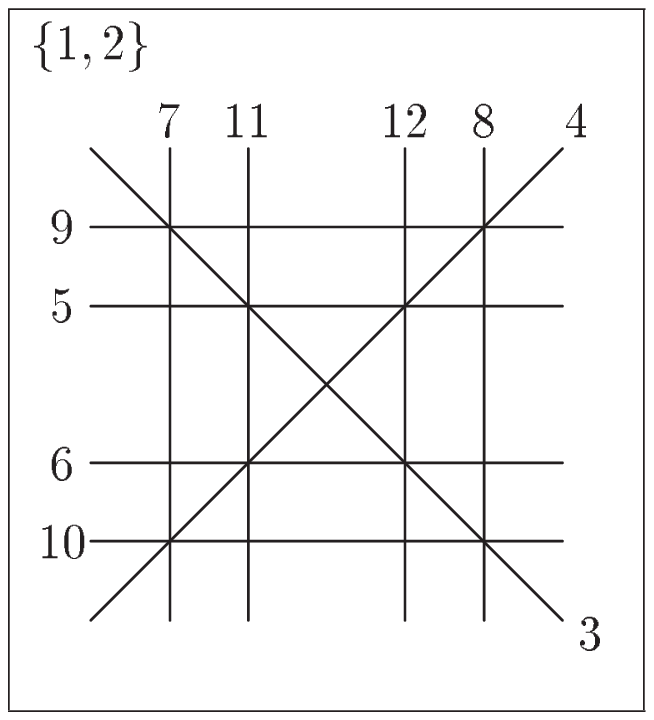} %
\qquad 
	\includegraphics[width=.4\textwidth]{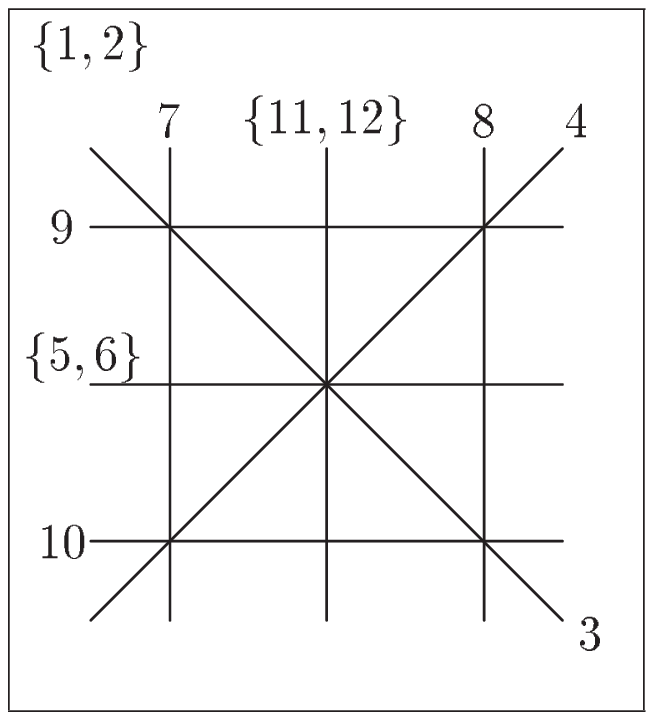} %
\caption{Degenerations  of Steiner's arrangement} 
\label{fig:dSteiner}
\end{figure}

The degeneration of $M[K_2]$ such that $1$ and $2$ are parallel, that is, 
$\{ 1, 2 \}$ is a circuit,  
has a realization defined by the left one in Figure \ref{fig:dSteiner}.   
Moreover, 
the degeneration of $M[K_2]$ such that 
$\{ 1, 2 \}$, $\{ 5, 6 \}$ and $\{ 11, 12 \}$ are circuits,  
is realizable. 
This realization is the $B_3$-arrangement (the right one in Figure \ref{fig:dSteiner}). 
Therefore, these two arrangements have 
weights with non-vanishing first cohomology 
in the same way of Remark \ref{rem:3}. 

\subsection{$m=3$ and $s=2$ (Libgober \cite{Li})} 
Two Latin squares 
$$  
K_1 = 
\begin{pmatrix} 
 1 & 2 & 3 \\ 
 3 & 1 & 2 \\ 
 2 & 3 & 1 
\end{pmatrix}, \qquad \text{ and } \qquad  
K_2 = 
\begin{pmatrix} 
 1 & 2 & 3 \\ 
 2 & 3 & 1 \\ 
 3 & 1 & 2 
\end{pmatrix}   
$$ 
are mutually orthogonal.  
We have 
$$ 
\C[K_1] = 
\begin{bmatrix}
  (1,4,7) & (1,5,8) & (1,6,9)  \\ 
  (2,4,9) & (2,5,7) & (2,6,8)  \\ 
  (3,4,8) & (3,5,9) & (3,6,7)  
\end{bmatrix}, \ 
\C[K_2] = 
\begin{bmatrix}
  (1,4,10) & (1,5,11) & (1,6,12) \\ 
  (2,4,11) & (2,5,12) & (2,6,10) \\ 
  (3,4,12) & (3,5,10) & (3,6,11)  
\end{bmatrix}.    $$ 
The matroid $M[K_1,K_2]$ is $\rm{AG}(2,3)$ (see \cite{Ox}) 
and realized as the Hessian configuration. 
The Hessian configuration is 
the arrangement of 12 projective lines 
passing through the nine inflection points 
of a nonsingular cubic 
in $\mathbb{P}^2(\mathbb{C})$ \cite[Example 6.30]{OT}, 
which we can define by lines 
$$ H_1= \{ x=0 \}, H_2= \{ y=0 \}, H_3= \{ z=0 \}, $$  
$$ H_4= \{ x+y+z=0 \}, 
   H_5= \{ x+ \omega^2 y + \omega z =0 \}, 
   H_6= \{ x+ \omega y + \omega^2 z =0 \}, $$  
$$ H_7= \{ x+ \omega y + \omega z =0 \}, 
   H_8= \{ x+   y + \omega^2 z =0 \}, 
   H_9= \{ x+ \omega^2 y +   z =0 \}, $$  
$$ H_{10}= \{ x+ \omega^2 y + \omega^2 z =0 \}, 
   H_{11}= \{ x+  \omega  y + z =0 \}, 
   H_{12}= \{ x+  y + \omega z =0 \}, $$  
where $\omega = e^{2 \pi i /3}$. 
Underlying matroids of arrangements 
$$ \{ H_1, \ldots, H_6, H_7, H_8, H_9 \}  
\ \text{ and } \  \{ H_1, \ldots, H_6, H_{10}, H_{11}, H_{12} \}  $$ 
are $M[K_1]$ and $M[K_2]$, respectively. 
The Hessian configuration $\{ H_1, \ldots, H_{12} \}$ 
has the underlying matroid $M[K_1,K_2]$ and  
we have 
$\dim{H^1(A(M[K_1, K_2]), e_{\lambda})} = 2$ 
for a non-zero one-form 
$$ e_{\lambda} = 
   \lambda_1 (e_1 + e_2 + e_3) 
 + \lambda_2 (e_4 + e_5 + e_6 ) 
 + \lambda_3 (e_7 + e_8 + e_9 ) 
 + \lambda_4 (e_{10} + e_{11} + e_{12} )   $$ 
with $\sum_{j=1}^{4} \lambda_j = 0$.

\subsection{Monomial arrangements (Cohen and Suciu \cite{CS})} 
Let $K$ be the Latin square of order $m$ 
defined by the addition table for $\mathbb{Z}_m \times \mathbb{Z}_m$ 
for $m \geq 2$.  
The monomial arrangement $\A_{m,m,3}$ in $\mathbb{C}^3$ is 
given by the defining polynomial 
$$  Q(\A_{m,m,3})  = (x_1^m - x_2^m) (x_1^m - x_3^m) (x_2^m - x_3^m).   $$ 
Set $\zeta = \exp{(2 \pi i/m)}$.  
Define  
$$ \A_{ij} = \{  H_{i,j}^{k} = \Ker{(x_i - \zeta^k x_j)}  : 1 \leq k \leq m  \}  $$ 
for $1 \leq i < j \leq 3$.   
So we have  $\A_{m,m,3} = \A_{12} \cup \A_{23} \cup \A_{13}$. 
Since $\cap_{k=1}^m H_{i,j}$ has rank two, 
the underlying matroid $M(\A_{ij})$ of $\A_{ij}$ is isomorphic to 
the uniform matroid $U_{2,m}$ of rank two.   
Another rank two intersections 
are $H_{1,2}^{p} \cap H_{2,3}^{q} \cap H_{1,3}^{r}$ for $p+q \equiv r \mod{m}$. 
Hence,  $K$ can be consider as the Latin square 
with rows indexed by $\A_{12}$, columns by $\A_{23}$, 
and symbols by $\A_{13}$ .  
The underlying matroid of $\A_{m,m,3}$  
is the matroid $M[K; M(\A_{12}),M(\A_{23}),M(\A_{13})]$.  
By Proposition \ref{prop:deg1}, 
$\A_{m,m,3}$  has weights with non-vanishing first cohomology. 

\subsection{Higher case ($\ell =3$)} 

Let $K$ be a Latin $3$-dimensional hypercube on $[2]$ defined by Figure \ref{3hypercube}. 
\begin{figure}[h]
\setlength{\unitlength}{0.45mm} 
\begin{picture}(130, 60)(0,0) 
 \put(0,30){ \framebox(20,10){ 1 } } 
 \put(0,0){ \framebox(20,10){ 2 } } 
 \put(40,0){ \framebox(20,10){ 1 } } 
 \put(40,30){ \framebox(20,10){ 2 } } 
 \put(70,47){ \framebox(20,10){ 2 } } 
 \put(70,17){ \framebox(20,10){ 1 } } 
 \put(110,17){ \framebox(20,10){ 2 } } 
 \put(110,47){ \framebox(20,10){ 1 } } 
     \put(22,5){\line(1,0){20}} %
     \put(22,35){\vector(1,0){20}} %
     \put(12,30){\vector(0,-1){20}} %
     \put(52,30){\line(0,-1){20}} %
     \put(92,22){\line(1,0){20}} %
     \put(92,52){\line(1,0){20}} %
     \put(80,47){\line(0,-1){20}} %
     \put(120,47){\line(0,-1){20}} %
    \put(22,5){\line(3,1){50}} %
    \put(22,35){\vector(3,1){50}} %
    \put(62,5){\line(3,1){50}} %
    \put(62,35){\line(3,1){50}} %
\end{picture} 
\hfill 
\begin{picture}(130, 60)(0,0) 
 \put(0,30){ \framebox(20,10){ 1357 } } 
 \put(0,0){ \framebox(20,10){ 2358 } } 
 \put(40,0){ \framebox(20,10){ 2457 } } 
 \put(40,30){ \framebox(20,10){ 1458 } } 
 \put(70,47){ \framebox(20,10){ 1368 } } 
 \put(70,17){ \framebox(20,10){ 2367 } } 
 \put(110,17){ \framebox(20,10){ 2468 } } 
 \put(110,47){ \framebox(20,10){ 1467 } } 
     \put(22,5){\line(1,0){20}} %
     \put(22,35){\vector(1,0){20}} %
     \put(12,30){\vector(0,-1){20}} %
     \put(52,30){\line(0,-1){20}} %
     \put(92,22){\line(1,0){20}} %
     \put(92,52){\line(1,0){20}} %
     \put(80,47){\line(0,-1){20}} %
     \put(120,47){\line(0,-1){20}} %
    \put(22,5){\line(3,1){50}} %
    \put(22,35){\vector(3,1){50}} %
    \put(62,5){\line(3,1){50}} %
    \put(62,35){\line(3,1){50}} %
\end{picture} 
	\caption{$K$ and $\C[K]$}   
	\label{3hypercube}
\end{figure} 
The matroid $M[K]$ is the matroid of type $L_8$ in \cite[{p.510}]{Ox}.  
Let $\A$  be an $4$-arrangement defined by the defining polynomial 
$$ x_1 x_2 x_3 x_4 (x_1+ x_2 + x_3 + x_4) (x_1+ bc x_2 + bx_3 + cx_4) 
           (x_1+ c x_2 + x_3 + c x_4) (x_1+ b x_2 + b x_3 + x_4),     $$ 
where $0, 1, b, c, bc$ are distinct each other. 
By the simple computation, 
$\A$ is a realization of $M[K]$. 
Therefore, $\A$ has weights with non-vanishing second cohomology. 
Let $\B$  be an $4$-arrangement defined by the defining polynomial 
$$ (x_1 - x_2)  (x_1 + x_2)  (x_2 - x_3) (x_2 + x_3) (x_3 - x_4) (x_3 + x_4) 
             (x_4 - x_1)  (x_4 + x_1).    $$ 
By the simple computation, 
we can check that $\B$ has no $3$-circuits and 
the family of $4$-circuits is 
$$  \C[K] \cup \{ (1,2,3,4), (1,2,7,8), (3,4,5,6), (5,6,7,8) \}. $$ 
Therefore, $\B$ has weights with non-vanishing second cohomology.

\begin{acknowledgments} 
A main part of this work was done while the author 
was visiting Mathematical Sciences Research Institute in August-November, 2004. 
The author would like to thank MSRI for its hospitality. 
The author is grateful to 
Professor Sergey Yuzvinsky and Professor Michael Falk 
for many helpful suggestions and valuable discussions and 
the author would also like to thank Professor Hiroaki Terao for his help. 
\end{acknowledgments}

\end{document}